\documentclass[journal]{IEEEtran}
\ifCLASSINFOpdf
\else
\fi
\usepackage[pdftex]{graphicx}
\usepackage[cmex10]{amsmath}
\usepackage[pdftex]{hyperref}

\usepackage[OT2,OT1]{fontenc}
\newcommand\cyr{%
\renewcommand\rmdefault{wncyr}%
\renewcommand\sfdefault{wncyss}%
\renewcommand\encodingdefault{OT2}%
\normalfont
\selectfont}
\DeclareTextFontCommand{\textcyr}{\cyr}

\begin{document}
\title{Nonlinear Control Synthesis\\ for a Self-energizing Electro-Hydraulic Brake}
\author{Alexey~Starykh,~\IEEEmembership{Tomsk~Polytechnic~University}
\thanks{A. Starykh is with the Department for Electric Drive and Electrical Equipment of Tomsk Polytechnic University, Lenin Ave. 30, 634050 Tomsk, Russia. E-mail: \texttt{starykh.alexey@sibmail.com}}}
\markboth{Journal "Izv.~Vuzov.~Electromechanics",~No.~6, December~2008}%
{Shell \MakeLowercase{\textit{et al.}}}
\maketitle

\begin{abstract}
Nonlinear control algorithm for a self-energizing electro-hydraulic brake is analytically designed. The desired closed-loop system behavior is reached via a synthesized nonlinear controller.
\end{abstract}


\IEEEpeerreviewmaketitle

\section{Introduction}
\IEEEPARstart{S}{elf}-reinforcing brakes are a subject of intensive investigation during last years. Working principle of such brakes is to use the wheelset's inertia momentum of a vehicle as the source of power for braking. One major advantage of self-reinforcing brakes is the energy consumption decrease which makes this research direction perspective.

At present the development of a new brake concept of a self-energizing electro-hydraulic brake for a railway application is being carried out at the Institute for Fluid Power Drives and Control (IFAS, RWTH Aachen University). Working principle of the braking system can be found in [6--9].

\section{Problem Statement}
The control task of the self-energizing electro-hydraulic brake is to track the reference signal of the pressure \(p_{sup}\) (the output variable) in the supply line of the brake system. The control system with a pure proportional controller allows reaching the goal, as shown in Fig. 1, where the reference signal represents a step function, whose values are the sequence \{27, 59, 91, 59, 27 bar\}, which corresponds to brake forces of \{5, 10, 15, 10, 5 kN\}. However, as seen from the figure, the oscillations of the supply line pressure appear in the system, which cause the oscillations of the brake force.

\begin{figure}[h!]
\center
\includegraphics{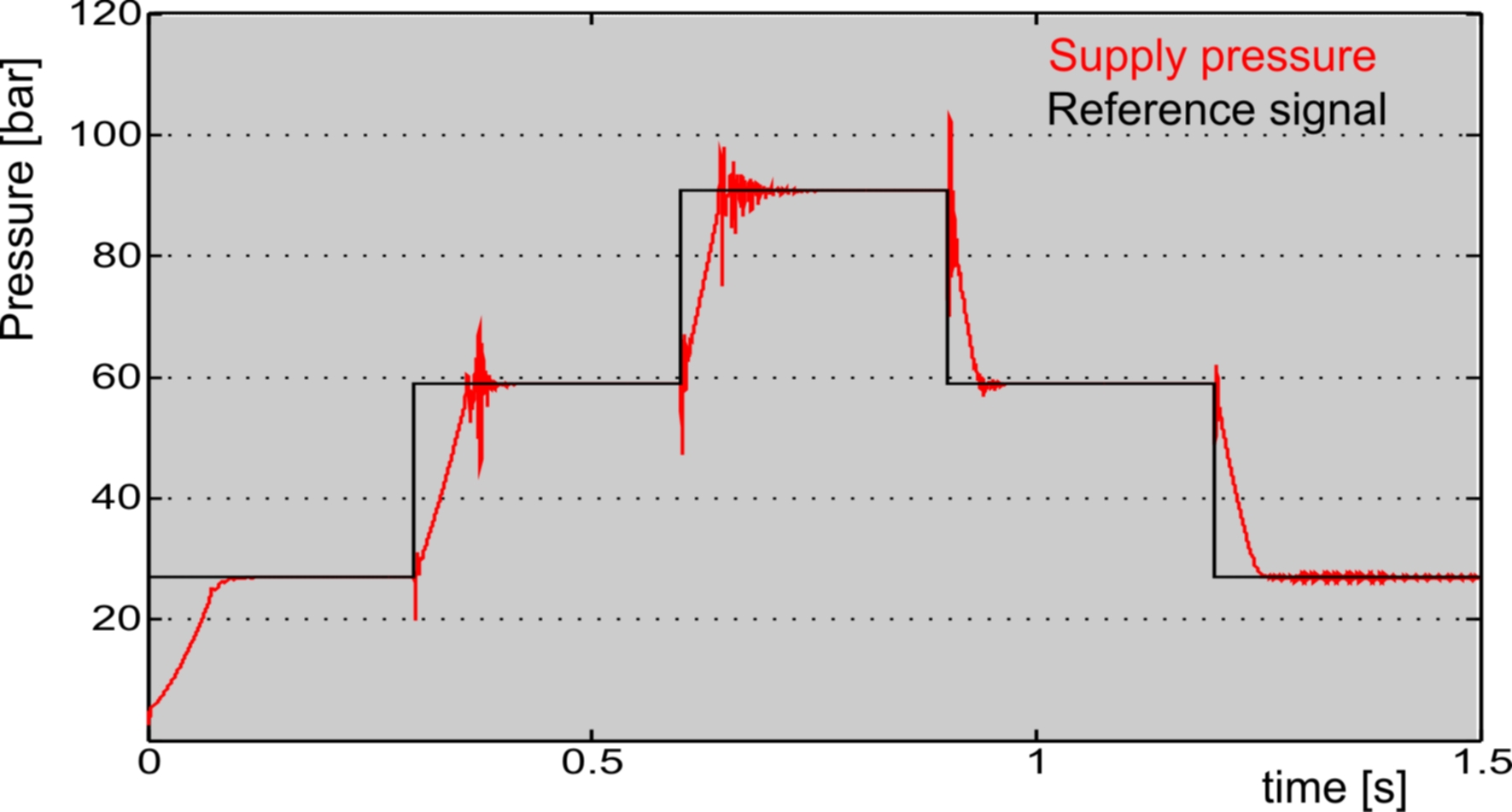}
\caption{The supply line pressure in the control system with a pure proportional controller}
\end{figure}

Such a behavior of braking, of course, is inadmissible for passenger trains. Therefore, the main requirement for the closed-loop system behavior is to get an aperiodic process of supply line pressure changes independent of the desired brake force level and friction coefficient variations of brake pads.

\section{Coordinates Transformation and Input-Output Linearization of a Plant Model}
To reach the required behavior of the controlled variable a controller of nonlinear structure is proposed to synthesize on basis of the exact feedback linearization method [3], [5].

Mathematical description of the electro-hydraulic system can be represented by a nonlinear model of 10\textsuperscript{th} order [6]. High order of the model result in the analytically complicated controller structure, which is a problem for simulation and experiments. Due to simplifications based on the results of the brake system analysis [6] as well as the approach used in [1], the basic model can be reduced to a 4\textsuperscript{th} order nonlinear model. For positive (\textbf{Case A}) and negative (\textbf{Case B}) valve opening this model has the form

\textbf{Case A:}
\noindent
\begin{equation}
\begin{pmatrix}
   \dot{p}_{L} \\
   \dot{p}_{sup} \\
   \dot{v}_{\textsf{v}} \\
   \dot{x}_{\textsf{v}} \\
\end{pmatrix}
   \!=\!
\begin{pmatrix}
   T^{A}_{L} x_{\textsf{v}} \sqrt{p_{sup} - p_{L} - \alpha p_{lp}} \\
   T^{A}_{sup} x_{\textsf{v}} \sqrt{p_{sup} - p_{L} - \alpha p_{lp}} \\
   -2 D_{\textsf{v}} \omega_{\textsf{v}} v_{\textsf{v}} - \omega^{2}_{\textsf{v}} x_{\textsf{v}} \\
   v_{\textsf{v}} \\
\end{pmatrix}
   \!+\!
\begin{pmatrix}
   0 \\
   0 \\
   \omega^{2}_{\textsf{v}} K_{\textsf{v}} \\
   0 \\
\end{pmatrix}
   \!u
\end{equation}

\textbf{Case B:}
\begin{equation}
\begin{pmatrix}
   \dot{p}_{L} \\
   \dot{p}_{sup} \\
   \dot{v}_{\textsf{v}} \\
   \dot{x}_{\textsf{v}} \\
\end{pmatrix}
   \!=\!
\begin{pmatrix}
   T^{B}_{L} x_{\textsf{v}} \sqrt{p_{L} + \alpha p_{sup} - p_{lp}} \\
   T^{B}_{sup} x_{\textsf{v}} \sqrt{p_{L} + \alpha p_{sup} - p_{lp}} \\
   -2 D_{\textsf{v}} \omega_{\textsf{v}} v_{\textsf{v}} - \omega^{2}_{\textsf{v}} x_{\textsf{v}} \\
   v_{\textsf{v}} \\
\end{pmatrix}
   \!+\!
\begin{pmatrix}
   0 \\
   0 \\
   \omega^{2}_{\textsf{v}} K_{\textsf{v}} \\
   0 \\
\end{pmatrix}
   \!u,
\end{equation}

\noindent where \(p_{L}, p_{sup}, v_{\textsf{v}}, x_{\textsf{v}}\) -- the load pressure of the brake actuator, the pressure in the supply line, the velocity and the control valve spool movement, respectively; \(\alpha\) -- the ratio between piston areas of the brake actuator; \(p_{lp}\) -- the pressure in the low pressure line; \(K_{\textsf{v}}, D_{\textsf{v}}, \omega_{\textsf{v}}\) -- the control valve parameters; \(T^{A}_{L}, T^{B}_{L}\) -- the known constant parameters; \(T^{A}_{sup}, T^{B}_{sup}\) -- the parameters dependent of the brake pads friction coefficient; \(u\) -- the control valve input signal. The output of the brake system is
\begin{equation}
   y=p_{sup}
\end{equation}
\noindent This model will be used for the controller synthesis in the following.

The nonlinear model (1) -- (3) can be represented as
\begin{flalign}
   \dot{x} & =f(x)+g(x)u\\
   y & =h(x)
\end{flalign}

The idea of the exact feedback linearization method is to find a nonlinear transformation (linearization algorithm) of a control signal, for which the model (4), (5) is linear or equivalent to a linear model in new coordinates \(z=H(x)\), where \(H(x)\) is the coordinates transformation.

An important notion of an input-output model is its relative degree. For linear systems, as known, the relative degree is the difference between the number of poles and zeros of the transfer function of a system. This is also the number of times an output needs to be differentiated in order that an input appears in the equation. For nonlinear systems it is defined in the similar manner. By differentiating the output \(y\) and substituting (4) we get (see [2])
\begin{equation}
\begin{split}
   y^{(1)}=\frac{\partial h}{\partial x}\frac{\partial x}{\partial t} & =\frac{\partial h}{\partial x}(f(x)+g(x)u) = {} \\
   & =L_{f}h(x)+L_{g}h(x)u,
\end{split}
\end{equation}
where \(L_{f}h\) and \(L_{g}h\) are Lie derivatives of the function \(h\) along the vector field \(f\) and \(g\), respectively.

For the electro-hydraulic brake, for positive and negative valve opening,
\(
   \frac{\partial h}{\partial x}g(x)=L_{g}h(x)=0
\)
for all \(x\) in the range of operating points. Therefore
\begin{equation*}
   y^{(1)}=L_{f}h(x)
\end{equation*}

Continuing in this way, we get
\begin{equation}
\begin{split}
   y^{(2)}=\frac{\partial L_{f}h}{\partial x}\frac{\partial x}{\partial t} & =\frac{\partial L_{f}h}{\partial x}(f(x)+g(x)u) = {} \\
   & =L^{2}_{f}h(x)+L_{g}L_{f}h(x)u,
\end{split}
\end{equation}
where \(L_{g}L_{f}h(x)=0\). Consequently,
\begin{equation*}
   y^{(2)}=L^{2}_{f}h(x)
\end{equation*}

The time derivative \(y^{(3)}\) yields:
\begin{equation}
   y^{(3)}=L^{3}_{f}h(x)+L_{g}L^{2}_{f}h(x)u,
\end{equation}
where the Lie derivative \(L_{g}L^{2}_{f}h(x)\ne0\) for all \(x\) in the range of operating points and for positive and negative valve opening has the form

\textbf{Case A:}
\begin{equation}
   L_{g}L^{2}_{f}h(x)=T^{A}_{sup} \sqrt{p_{sup} - p_{L} - \alpha p_{lp}} \: \omega^{2}_{\textsf{v}} K_{\textsf{v}}
\end{equation}

\textbf{Case B:}
\begin{equation}
   L_{g}L^{2}_{f}h(x)=T^{B}_{sup} \sqrt{p_{L} + \alpha p_{sup} - p_{lp}} \: \omega^{2}_{\textsf{v}} K_{\textsf{v}}
\end{equation}

Thus, the equation (8) with the nonzero factor for \(u\) describes the relation between the input \(u\) and the output \(y\). Here, according to the definition in [3], the relative degree of the 4th order nonlinear model of the electro-hydraulic brake~is~3.

For linear state-space systems, derivatives of an output are chosen as state variables of a plant. Using the similar approach, we determine new coordinates as
\begin{equation}
   z=
\begin{pmatrix}
   z_{1} \\
   z_{2} \\
   z_{3} \\
\end{pmatrix}
   =
\begin{pmatrix}
   y \\
   y^{(1)} \\
   y^{(2)} \\
\end{pmatrix}
   =
\begin{pmatrix}
   h(x) \\
   L_{f}h(x) \\
   L^{2}_{f}h(x) \\
\end{pmatrix}
   \doteq H(x)
\end{equation}

For positive and negative valve spool movement, the new coordinates \(z\) can be expressed by the coordinates \(x\) of the simplified brake model (1) -- (3), i.e.

\textbf{Case A:}
\begin{flalign*}
   z_{1} & =p_{sup} \\
   z_{2} & =T^{A}_{sup} x_{\textsf{v}} \sqrt{p_{sup} - p_{L} - \alpha p_{lp}} \\
   z_{3} & =\frac{1}{2} T^{A}_{sup} (T^{A}_{sup} - T^{A}_{L}) x^{2}_{\textsf{v}} + T^{A}_{sup} v_{\textsf{v}} \sqrt{p_{sup} - p_{L} - \alpha p_{lp}}
\end{flalign*}

\textbf{Case B:}
\begin{flalign*}
   z_{1} & =p_{sup} \\
   z_{2} & =T^{B}_{sup} x_{\textsf{v}} \sqrt{p_{L} + \alpha p_{sup} - p_{lp}} \\
   z_{3} & =\frac{1}{2} T^{B}_{sup} (\alpha T^{B}_{sup} - T^{B}_{L}) x^{2}_{\textsf{v}} + T^{A}_{sup} v_{\textsf{v}} \sqrt{p_{L} + \alpha p_{sup} - p_{lp}}
\end{flalign*}

In accordance with (11), equations (5) -- (7) can be rewritten in the form
\begin{equation}
   \dot{z}_{1}=z_{2},\quad \dot{z}_{2}=z_{3},\quad \dot{z}_{3}=a(x)+b(x)u,\quad y=z_{1},
\end{equation}
where
\begin{equation*}
   a(x)=L^{3}_{f}h(x),\quad b(x)=L_{g}L^{2}_{f}h(x)
\end{equation*}
The function \(b(x)\), for positive and negative valve spool movement, is determined by equations (9), (10), the function \(a(x)\) is expressed as

\textbf{Case A:}
\begin{equation*}
\begin{split}
   a(x) & =\frac{3}{2} T^{A}_{sup} (T^{A}_{sup} - T^{A}_{L}) v_{\textsf{v}} x_{\textsf{v}} - {} \\
& - T^{A}_{sup} (2 D_{\textsf{v}} \omega_{\textsf{v}} v_{\textsf{v}} + \omega^{2}_{\textsf{v}} x_{\textsf{v}}) \sqrt{p_{sup} - p_{L} - \alpha p_{lp}}
\end{split}
\end{equation*}

\textbf{Case B:}
\begin{equation*}
\begin{split}
   a(x) & =\frac{3}{2} T^{B}_{sup} (\alpha T^{B}_{sup} + T^{B}_{L}) v_{\textsf{v}} x_{\textsf{v}} - {} \\
& - T^{B}_{sup} (2 D_{\textsf{v}} \omega_{\textsf{v}} v_{\textsf{v}} + \omega^{2}_{\textsf{v}} x_{\textsf{v}}) \sqrt{p_{L} + \alpha p_{sup} - p_{lp}}
\end{split}
\end{equation*}

If models (1) -- (3), and (12) where equivalent, than the exact feedback linearization problem would be solvable. This, in its turn, would mean that a control signal of the form
\begin{equation}
   u=\frac{1}{b(x)}(-a(x)+\nu),
\end{equation}
where \(\nu\) is a new control signal, due to overall compensation of the nonlinear functions \(a(x)\) and \(b(x)\), would lead the plant model (1) -- (3), on the assumption of parametric certainty of the model and measurability of the state variables \(x\), to a system, whose behaviour is exactly identical with the behaviour of the linear model
\begin{flalign}
   \dot{z}= 
\begin{pmatrix}
   0 & 1 & 0 \\
   0 & 0 & 1 \\
   0 & 0 & 0 \\
\end{pmatrix}
\begin{pmatrix}
   z_{1} \\
   z_{2} \\
   z_{3} \\
\end{pmatrix}&+
\begin{pmatrix}
   0 \\
   0 \\
   1 \\
\end{pmatrix}\!\nu, \\
   y=
\begin{pmatrix}
   1 & 0 & 0 \\
\end{pmatrix}&
\begin{pmatrix}
   z_{1} \\
   z_{2} \\
   z_{3} \\
\end{pmatrix}\!,
\end{flalign}
whose output \(y=z_{1}\) would coincide with the output \(y=p_{sup}\) of the simplified nonlinear model.

However, since rank~\(H(x)=3\) for all \(x\) in the range of operating points, which is less then the dimension \(n=4\) of the model (1) -- (3), the mapping \(H(x)\) is not a diffeomorphism. Let us note that a mapping \(\varphi \colon X\to Y\), where \(X\) and \(Y\) are smooth manifolds of dimension \(n\), is a diffeomorphism if \(\varphi\) is bijective and both \(\varphi\) and \(\varphi^{-1}\) are smooth mapping (see Appendix A in [3]).

Considering aforesaid, the model (12) is not equivalent to the simplified nonlinear model (1) -- (3). This is caused by the internal dynamics, i.e. behaviour of \((n-3)\)-dimensional part of the nonlinear model, which, however, has no influence on the output (15) of the linear model (14).

In that case the problem of input-output linearization can be solved [4]. This means that the control (13) transforms the nonlinear model (1) -- (3) of the electro-hydraulic brake into a model whose input-output behaviour can be represented in the form (14), (15). Additionally, one more variable \(\eta=T(x)\) has to be supplemented with the new coordinates \(z\). According to Proposition 4.1.3 in [3] it is always possible to find such a function that the jacobian matrix of the mapping
\begin{equation*}
   M(x)=
\begin{pmatrix}
   T(x) \\
   H(x) \\
\end{pmatrix}
\end{equation*}
\noindent is nonsingular at some point \(x=x^{*}\).

For the plant in question the load pressure of the brake actuator \(p_{L}\) was chosen as the additional variable, i.e. \(\eta=T(x)=p_{L}\). The jacobian matrix of the mapping \(M(x)\) is nonsingular at each operating point of the brake. The time derivative of \(\eta\) for positive and negative valve opening has the form

\textbf{Case A:}
\begin{equation}
   \dot{\eta}\doteq q(z)=\frac{T^{A}_{L}}{T^{A}_{sup}}z_{2}
\end{equation}

\textbf{Case B:}
\begin{equation}
   \dot{\eta}\doteq q(z)=\frac{T^{B}_{L}}{T^{B}_{sup}}z_{2}
\end{equation}
Thus, the equivalent model in the new coordinates can be written as
\begin{flalign}
   \dot{\eta}&=q(z),\\
   \dot{z}=
\begin{pmatrix}
   0 & 1 & 0 \\
   0 & 0 & 1 \\
   0 & 0 & 0 \\
\end{pmatrix}
\begin{pmatrix}
   z_{1} \\
   z_{2} \\
   z_{3} \\
\end{pmatrix}
&+
\begin{pmatrix}
   0 \\
   0 \\
   1 \\
\end{pmatrix}
   (\overline{a}(\eta,z)+\overline{b}(\eta,z)u),\nonumber\\
   y&=z_{1},\nonumber
\end{flalign}
where the equation (18) for \mbox{\(z\equiv0\)} describes the zero-dynamics of the model (1) -- (3). Taking into consideration that \(x=M^{-1}(\eta,z)\), the functions \(\overline{a}(\eta,z)\) and \(\overline{b}(\eta,z)\) are expressed as
\begin{equation*}
   \overline{a}(\eta,z)=a\{M^{-1}(\eta,z)\},\quad \overline{b}(\eta,z)=b\{M^{-1}(\eta,z)\}
\end{equation*}
\noindent The block diagram of the equivalent model is depicted in Fig.~2.

\begin{figure}[h!]
\center
\includegraphics{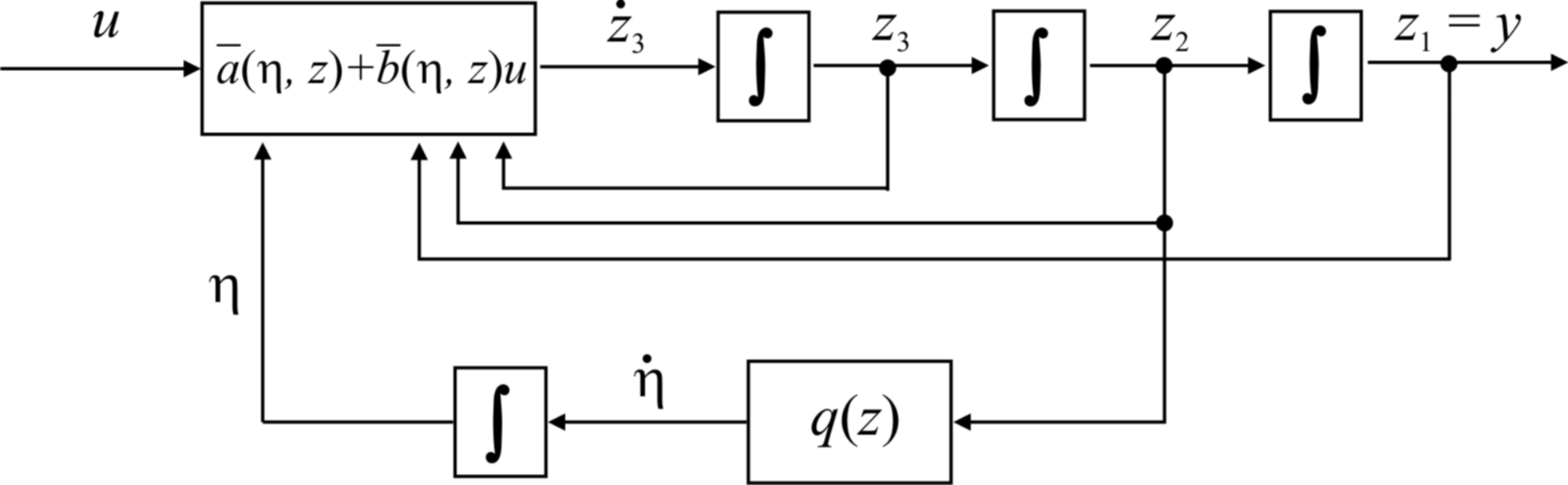}
\caption{The normal form of the electro-hydraulic brake model (1) -- (3)}
\end{figure}

Thus we can conclude that for the self-energizing electro-hydraulic brake the controller synthesis based on the exact feedback linearization method is conventionally divided into the following steps. Firstly, the control signal \(\nu\) is determined for the linear model (14), (15), which describes the relationship between input and output variables of the nonlinear model (1) -- (3). Secondly, we take into consideration the internal dynamics of the plant, described by the equation (18), which has no influence on the output of the linear model, see Fig.~3.

\begin{figure}[h!]
\center
\includegraphics{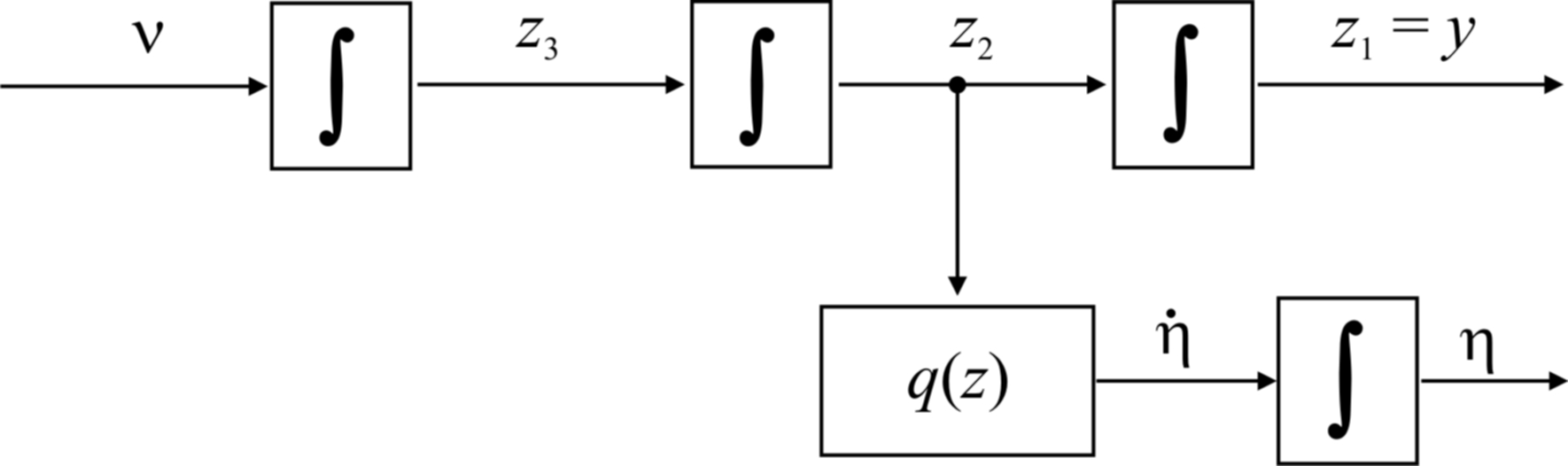}
\caption{The linearized model and the internal dynamics}
\end{figure}

Finally, we generate the control signal \(u\) of the 4\textsuperscript{th} order nonlinear model using the nonlinear transformation for the signal \(\nu\). The control \(u\) is expressed as
\begin{equation*}
   u=\frac{1}{\overline{b}(\eta,z)}(-\overline{a}(\eta,z)+\nu)
\end{equation*}

The obtained signal \(u\) is the input signal for the real brake system, described by the nonlinear model of 10\textsuperscript{th} order.

\section{Control Synthesis for the Linear Model of the Electro-Hydraulic Brake}
The linear model (14) gives the opportunity to use classical methods of linear systems theory to achieve the required control performance, i.e. to get an aperiodic behaviour of the output variable (15) of the model (14). One such approach is the pole placement method which allows reaching any prescribed placement of closed loop system poles [10].

According to the method the control input signal \(\nu\) for the linear model (14) is chosen in the form
\begin{equation*}
   \nu=K_{z}e,
\end{equation*}
where \mbox{\(e=z^{*}-z\)} is the error, i.e. the difference between the reference state vector and the real state vector. Since the reference signal \(p^{*}_{sup}\) is assumed to be a piecewise constant function, the reference signal \(z^{*}\) of the model (14) during each time period, where \(p^{*}_{sup}\) is constant, can be represented as
\begin{equation*}
   z^{*}=
\begin{pmatrix}
   p^{*}_{sup} \\
   \dot{p}^{*}_{sup} \\
   \ddot{p}^{*}_{sup} \\
\end{pmatrix}
   =
\begin{pmatrix}
   p^{*}_{sup} \\
   0 \\
   0 \\
\end{pmatrix}
\end{equation*}
The row-vector \(K_{z}\) provides prescribed poles placement of the following closed loop system
\begin{equation*}
   \dot{e}=(A - b K_{z})e,
\end{equation*}
where
\begin{equation*}
   A=
\begin{pmatrix}
   0 & 1 & 0 \\
   0 & 0 & 1 \\
   0 & 0 & 0 \\
\end{pmatrix}
;\quad b=
\begin{pmatrix}
   0 \\
   0 \\
   1 \\
\end{pmatrix}
\end{equation*}

Since the electro-hydraulic brake system in question is a SISO system, i.e. a single-input single-output system, only the output \(p_{sup}\) is measured. Consequently, not all state vector \(x\) is known. It means that both the state vector \(z\) of the model (14), (15) and the variable \(\eta\) of the model (18) are unmeasured. Instead of these variables their estimates have to be used. Therefore, the control system is supplemented with a full-order state observer for (14), (15). By means of the estimate \(\hat{z}\) of the vector \(z\) the estimate \(\hat{\eta}\) of the variable \(\eta\) can be calculated (see (16) and (17)). The full-order state observer model is determined on the basis of the model (14), (15) and has the form:
\begin{flalign*}
   \dot{\hat{z}}&=A\hat{z} + b\nu + K_{obs}(y-\hat{y}), \\
   \hat{y}&=c^{T}\hat{z},
\end{flalign*}
where \(K_{obs}\) is a column-vector of preset coefficients and \mbox{\(c^{T}=(1\;0\;0)\)}. The vector \(K_{obs}\) is chosen such that the matrix of closed loop estimation system \((A-K_{obs}c^{T})\) be Hurwitz stable [10].

Simulation results of the electro-hydraulic brake system of 10\textsuperscript{th} order with the controller of described nonlinear structure are shown in Fig.~4.

\begin{figure}[h!]
\center
\includegraphics{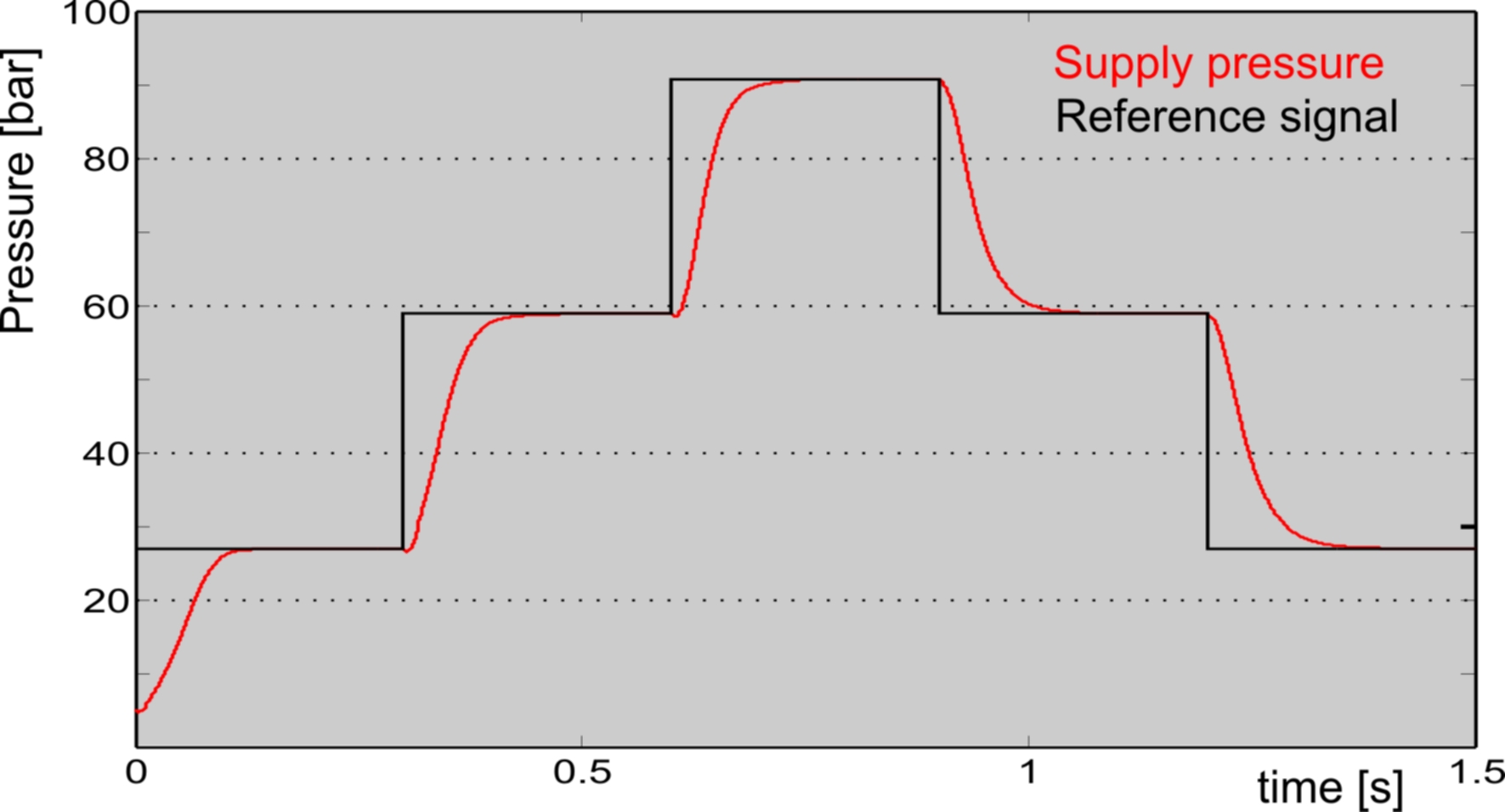}
\caption{The supply line pressure in the control system with the nonlinear controller}
\end{figure}

One can see from the figure, that the nonlinear control algorithm synthesized on the basis of the exact feedback linearization method allows to get the desired closed-loop system behavior. The results of accomplished simulation demonstrated as well, that the control system is robust with respect to small variations of the friction coefficient. Expansion of the robust stability ranges stimulates further development of the controller for the self-energizing electro-hydraulic brake.

\section{Conclusion}
Nonlinear control algorithm synthesis of the self-energizing electro-hydraulic brake on the basis of the exact feedback linearization method has been fulfilled. With the help of the obtained algorithm the problem of the desired closed-loop system behavior has been solved.

\section*{Acknowledgment}

The author would like to thank German Academic Exchange Service (DAAD) and Ministry of Education and Science of Russian Federation for supporting the research stay at the Institute for Fluid Power Drives and Control in Aachen (Project No. 7.375.C2007).

\ifCLASSOPTIONcaptionsoff
  \newpage
\fi

\end{document}